\documentclass[10pt,leqno]{amsart}
\usepackage{indentfirst, amssymb,amsmath,amsthm}
\usepackage{hyperref}
\usepackage{csquotes}
\usepackage{graphicx}
\usepackage{epstopdf}
\usepackage[margin=.75in]{geometry}
\usepackage{tikz,pgfplots}
\usepgfplotslibrary{fillbetween}

\newtheorem{cor}{Remark}
\newcommand{\beas}{\begin{eqnarray*}}
\newcommand{\eeas}{\end{eqnarray*}}
\begin{document}
\title[Geometric visualizations of $b^{e}<e^{b}$ when $e<b$]{Geometric visualizations of $b^{e}<e^{b}$ when $e<b$}
\date{}
\author[B. Chakraborty, S. Chakraborty and R. Farhadian]{Bikash Chakraborty, Sagar Chakraborty and Reza Farhadian}
\date{}
\address{Department of Mathematics, Ramakrishna Mission Vivekananda
Centenary College, Rahara, West Bengal 700 118, India.}
\email{bikashchakraborty.math@yahoo.com, bikash@rkmvccrahara.org}
\address{Department of Mathematics, Jadavpur University, Kolkata - 700032.}
\email{sagarchakraborty55@gmail.com}
\address{Department of Statistics, Lorestan University, Khorramabad, Iran}
\email{farhadian.reza@yahoo.com}
\maketitle
\begin{abstract}
In connection to the two fascinating constants $e$ and $\pi$, there are many beautiful visual proofs to the inequality $\pi^{e}<e^{\pi}$.  The aim of this classroom capsule is to give three visual proofs to the more general inequality $b^{e}<e^{b}$ where $e<b$.
\end{abstract}
\section{Introduction}
In connection to the two fascinating constants $e$ and $\pi$, there are many beautiful visual proof to the inequality $\pi^{e}<e^{\pi}$. The inequality $\pi^{e}<e^{\pi}$ has numerous geometric visualizations in the literature. In 1987, Nakhli (\cite{F}) gave a beautiful geometric visualization to the inequality $\pi^{e}<e^{\pi}$ using the fact that the global maximum of the curve $y=\frac{\ln x}{x}$ occurs at the point $x=e$. The same inequality was visualized by Nelsen (\cite{R}) using the fact that the curve $y=e^{\frac{x}{e}}$ lies above the line $y=x$. Recently, the first Author(\cite{C}) visualized the inequality $\pi^{e}<e^{\pi}$ from the Napier's inequality. Moreover, the first Author together with Mukherjee(\cite{MC}) gave another nice visualization of $\pi^{e}<e^{\pi}$ using the fact that the line $y=x-1$ is a tangent line to the curve $y=\ln x$ at the point $(1,0)$.\par
More recently, using the area argument, Haque (\cite{H}) gave a beautiful geometric visualization to the inequality $e^{b}> b^{e}$ when $e<b$. \par
In this classroom note, We provide three alternative geometric visualizations to the inequality $e^{b}> b^{e}$ where $e<b$.
\begin{center}
\section{Visual Proof 1}
Using the elementary  calculus, one can prove that $e^{x}>1+x$ for $x>0$ (see, \cite{BS}, pp. 191), and by putting $x=\frac{\pi}{e}$, one can obtain that  $\pi^{e}<e^{\pi}$. In this section, we will see the geometric visualization to the inequality $e^{x}>1+x$ for $x>0$, by  using the fact that the curve $y=e^{x}$ lies above the tangent line $y=x+1$ to the curve at the point $(0,1)$.
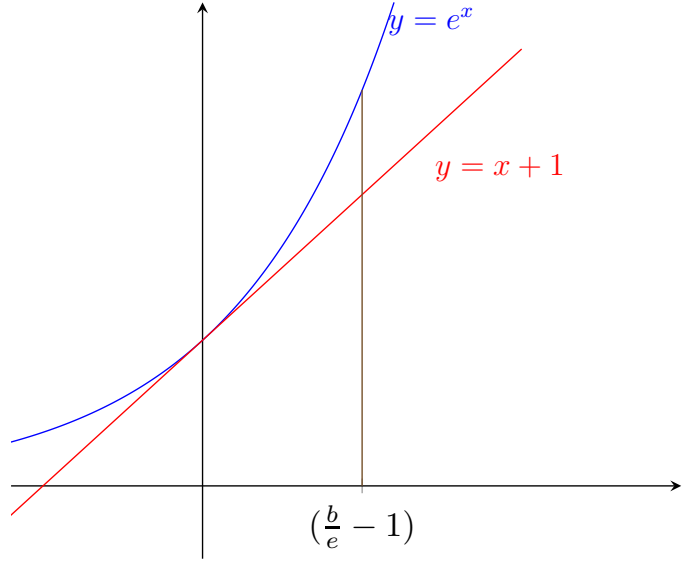
\begin{tikzpicture}[scale=1.3]
  \begin{axis}[
    axis lines=center,
    ymin=-.5,
    xmin=-1.2,
    xmax=3,
    no marks,
     xtick={0,1},
    xticklabels={,$(\frac{b}{e}-1)$},
    ytick=\empty,
    yticklabels={};
    ]
    \addplot+[smooth,blue, domain=-3:1.2, samples=1000] {e^x}; 
        \addplot+[smooth,red, domain=-2:2, samples=1000] {x+1}; 
\addplot +[mark=none] coordinates {(1, 0) (1, e)};
  \end{axis}
  \node[] at (4.3,5.5) {\Large\color{blue}$y=e^{x}$};
    \node[] at (5,4) {\Large\color{red}$y=x+1$};
\end{tikzpicture}
\end{center}
Since the line $x=\frac{b}{e}-1$ when $b>e$ intersects the line $y=x+1$ and the curve $y=e^{x}$ at $(\frac{b}{e}-1, \frac{b}{e})$ and $(\frac{b}{e}-1, e^{\frac{b}{e}-1})$ respectively, so,
$$\frac{b}{e}<e^{\frac{b}{e}-1}\Rightarrow b^{e}<e^{b}.$$
\section{Visual Proof 2}
Since the area of the rectangle bounded by $x$-axis, $y$-axis and the lines $y=1$,  $x=\ln (\frac{b}{e})$ is $1\cdot \ln (\frac{b}{e})$; and the area under the curve $y=e^{x}$ from $x=0$ to $x=\ln(\frac{b}{e})$ is $\int_{0}^{\ln(\frac{b}{e})}e^{x} dx$, thus from the corresponding figure, we can visualize that
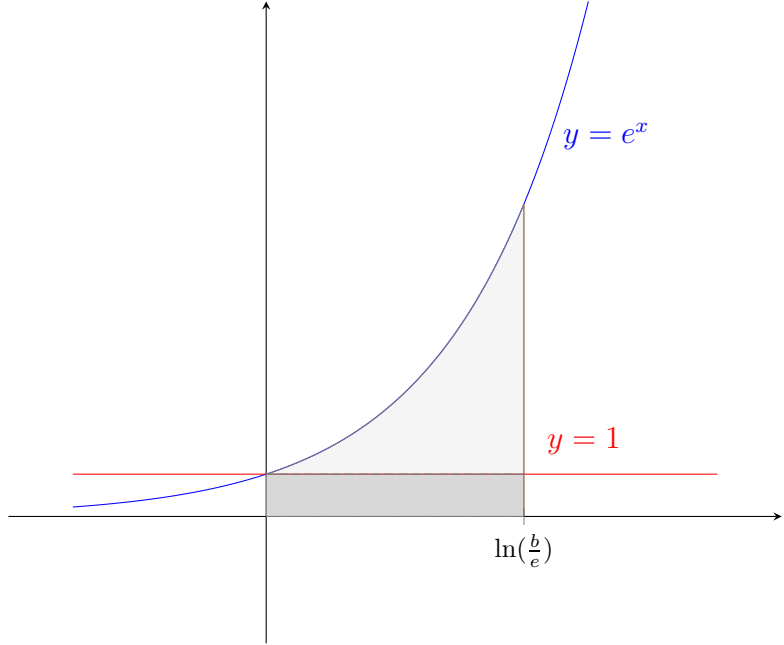
\begin{center}
\begin{tikzpicture}[scale=1.5]
  \begin{axis}[
    axis lines=center,
    ymin=-3,
    xmin=-2,
    xmax=4,
    no marks,
     xtick={0,2},
    xticklabels={,$\ln (\frac{b}{e})$},
    ytick=\empty,
    yticklabels={};
    ]
    \addplot+[smooth,blue, domain=-1.5:2.5, samples=1000] {exp(x)}; 
    \addplot+[smooth,red, domain=-1.5:3.5, samples=1000] {1}; 
\addplot +[mark=none] coordinates {(2, 0) (2, 7.38)};
      \addplot+[mark=none, gray, name path=A, domain=0:2, samples=500] {exp(x)}; 
      \addplot+[mark=none, gray, name path=B, domain=0:2, samples=500] {1}; 
    \addplot+[lightgray, opacity=.15] fill between[of=A and B]; 
      \addplot+[mark=none, gray, name path=C, domain=0:2, samples=500] {0}; 
      \addplot+[mark=none, gray, name path=B, domain=0:2, samples=500] {1}; 
    \addplot+[lightgray, opacity=.6] fill between[of=B and C]; 
  \end{axis}
  \node[] at (5.3,4.5) {\Large\color{blue}$y=e^{x}$};
  \node[] at (5.1,1.8) {\Large\color{red}$y=1$};
\end{tikzpicture}
\end{center}
$$\ln\left(\frac{b}{e}\right)<\int_{x=0}^{\ln(\frac{b}{e})}e^{x} dx=\frac{b}{e}-1\Rightarrow b^{e}<e^{b}.$$
\section{Visual Proof 3}
Since the area bounded by the lines $x=e$,  $x=b$, $y= \ln b$ and the curve $y=\ln x$ is $(b-e)\cdot\ln b -\int_{e}^{b} \ln x dx$, thus from the corresponding figure, we can visualize that
  
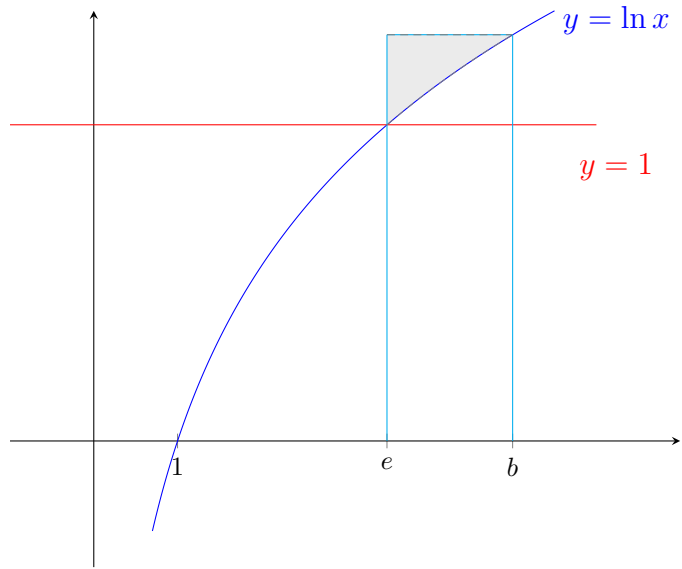
\begin{center}
\begin{tikzpicture}[scale=1.3]
  \begin{axis}[
    axis lines=center,
    ymin=-0.5,
    xmin=-1,
    xmax=7,
    no marks,
     xtick={0,1,3.5,5},
    xticklabels={,$1$,$e$,$b$},
    ytick=\empty,
    yticklabels={};
    ]
    \addplot+[smooth, blue, domain=0.7:5.5, samples=1000] {ln(x)}; 
    \addplot+[smooth, red, domain=-1:6, samples=1000] {1.253}; 
    \addplot+[smooth, cyan, domain=3.5:5, samples=1000] {1.61}; 
    \addplot +[mark=none, cyan] coordinates {(3.5, 0) (3.5, 1.61)};
    \addplot +[mark=none, cyan] coordinates {(5, 0) (5, 1.61)};
      \addplot+[mark=none, gray, name path=A, domain=3.5:5, samples=500] {1.61}; 
      \addplot+[mark=none, gray, name path=B, domain=3.5:5, samples=500] {ln(x)}; 
    \addplot+[lightgray, opacity=.3] fill between[of=A and B]; 
  \end{axis}
  \node[] at (6.2,5.6) {\Large\color{blue}$y=\ln x$};
\node[] at (6.2,4.1) {\Large\color{red}$y=1$};
\end{tikzpicture}
\end{center}
 \beas \left[(b-e)\cdot\ln b-\int_{e}^{b}\ln x ~dx \right]>0.\eeas
i.e.,\beas  b -e\ln b>0,\eeas
Thus \beas b^{e}<e^{b}.\eeas
\begin{center}
\begin{cor}
  The visual proof 3 can be made more general to prove the inequality $b^{a}$ and $a^{b}$, where $e\leq a<b$.
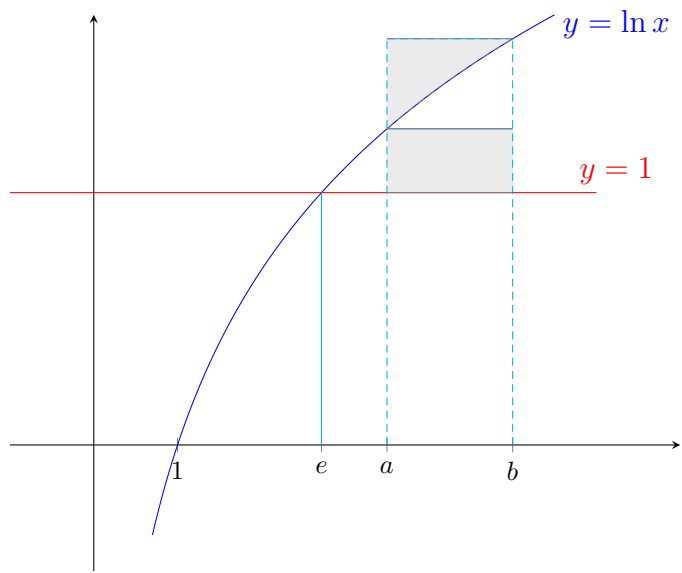
  \begin{center}
\begin{tikzpicture}[scale=1.3]
  \begin{axis}[
    axis lines=center,
    ymin=-0.5,
    xmin=-1,
    xmax=7,
    no marks,
     xtick={0,1,e,3.5,5},
    xticklabels={,$1$,$e$,$a$,$b$},
    ytick=\empty,
    yticklabels={};
    ]
    \addplot+[smooth, blue, domain=0.7:5.5, samples=1000] {ln(x)}; 
    \addplot+[smooth, red, domain=-1:6, samples=1000] {1}; 
    \addplot+[smooth, cyan, domain=3.5:5, samples=1000] {1.253}; 
    \addplot+[smooth, cyan, domain=3.5:5, samples=1000] {1.61}; 
    \addplot +[mark=none, cyan] coordinates {(e, 0) (e, 1)};
    \addplot +[mark=none, cyan] coordinates {(3.5, 0) (3.5, 1.61)};
    \addplot +[mark=none, cyan] coordinates {(5, 0) (5, 1.61)};
      \addplot+[mark=none, gray, name path=A, domain=3.5:5, samples=500] {1.61}; 
      \addplot+[mark=none, gray, name path=B, domain=3.5:5, samples=500] {ln(x)}; 
    \addplot+[lightgray, opacity=.3] fill between[of=A and B]; 
          \addplot+[mark=none, gray, name path=C, domain=3.5:5, samples=500] {1.253}; 
          \addplot+[mark=none, gray, name path=D, domain=3.5:5, samples=500] {1}; 
    \addplot+[lightgray, opacity=.3] fill between[of=C and D]; 
  \end{axis}
  \node[] at (6.2,5.6) {\Large\color{blue}$y=\ln x$};
\node[] at (6.2,4.1) {\Large\color{red}$y=1$};
\end{tikzpicture}
\end{center}
The shaded are described in the figure is \beas \left[(b-a)\cdot\ln b-\int_{a}^{b}\ln x ~dx \right]+\left[(b-a)\cdot\ln a -(b-a)\cdot 1\right],\eeas which is  strictly positive (It is seen from the figure). Since \beas \left[(b-a)\cdot\ln b-\int_{a}^{b}\ln x ~dx \right]+\left[(b-a)\cdot\ln a -(b-a)\cdot 1\right]&=& b\ln a -a\ln b,\eeas
Thus \beas b^{a}<a^{b}.\eeas
\end{cor}
\textbf{Disclosure statement}
\end{center}
No potential conflict of interest was reported by the authors.

\end{document}